\documentclass[a4paper]{article}

\usepackage[english]{babel}
\usepackage[utf8x]{inputenc}
\usepackage{amsmath,amsthm,amssymb}
\usepackage{etex}
\usepackage{graphicx}
\usepackage[colorinlistoftodos]{todonotes}
\usepackage{amsmath,pictexwd,verbatim}

\newtheorem{theorem}{Theorem}
\newtheorem{lemma}[theorem]{Lemma}

\title{A Direct Construction of Non-Transitive Dice Sets}
\author{Levi Angel and Matt Davis\footnote{mattd@muskingum.edu} \\ Muskingum University}
\begin{document}
\maketitle

\begin{abstract}
In this paper, we give a direct construction for a set of dice realizing any given tournament $T$. The construction for a tournament with $n$ vertices requires a number of sides on the order of $n$, which is the best general construction to date. Our construction relies only on a standard theorem from graph theory.
\end{abstract}

\section{Introduction}

Non-transitive dice have been a object of wide interest since Martin Gardner introduced some work (now known as Efron Dice) of Bradley Efron to the general public in \cite{G}.  For a set of at least 3 dice $\{X_{1},X_{2},X_{3}, \ldots\}$, with faces labeled in a nonstandard way, we define the relation $\succeq$ by declaring that $X_{i} \succeq X_{j}$ exactly if, when the dice are rolled, the probability that $X_{i}$ rolls a higher number than $X_{j}$ is greater than $1/2$. Paradoxically, it is entirely possible to put numbers on the dice in order to make the relation $\succeq$ non-transitive. The most basic example of such a set of dice is the following set of three 3-sided dice:

\begin{equation} \label{eq:example} \begin{tabular}{c|c} \textrm{Die:} & \textrm{Faces:} \\ \hline $X_{1}$ & 1,5,9 \\ $X_{2}$ & 3,4,8 \\ $X_{3}$ & 2,6,7 \end{tabular} \end{equation}

One can easily check that $X_{1} \succeq X_{2}$, $X_{2} \succeq X_{3}$, and $X_{3} \succeq X_{1}$, where the stronger die in each pair has probability 5/9 of winning. Note that we could widen our focus and examine non-transitive sets of general random variables. (See for example \cite{ST} and \cite{T}, which seem to predate any notion of non-transitive dice specifically). However, even focusing on dice, there are many open questions. %See \cite{C} for a recent example of work on non-transitive dice which is quite different from our work here.

\begin{figure}[ht]
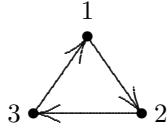

\[ \beginpicture
\setcoordinatesystem units <0.2in,0.2in>
\put{$\bullet$} at 0 5
\put{$\bullet$} at 1.4 3
\put{$\bullet$} at -1.4 3
\put{$1$} at 0 5.7
\put{$2$} at 1.9 3
\put{$3$} at -1.9 3
\arrow <10pt> [.2,.67] from 0 5 to 1.4 3
\arrow <10pt> [.2,.67] from 1.4 3 to -1.4 3
\arrow <10pt> [.2,.67] from -1.4 3 to 0 5
\endpicture \]
\caption{A non-transitive tournament on 3 vertices}
\label{tournex}
\end{figure}

Our aim in this paper is to consider larger sets of dice with arbitrary relations between them. A \textit{tournament} on $n$ vertices is a directed realization of the complete graph $K_{n}$. In other words, it is a directed graph on the vertices $\{1,2, \ldots n \}$, where for any pair of vertices $i$ and $j$, either there is an edge from $i$ to $j$ or from $j$ to $i$, but not both. We can interpret this as a definition of a relation on a set of dice - we say that a set of dice \textit{realizes} a tournament $T$ if $X_{i} \succeq X_{j}$ if and only if there is an edge from $i$ to $j$ in $T$. So the set of dice in \eqref{eq:example} realizes the tournament in Figure \ref{tournex}.

A natural problem to solve is, given an arbitrary tournament $T$, to construct a set of dice which realize $T$. There are many examples of general constructions that can be used to solve this problem. But going further, can we find a set of dice with a relatively small number of sides that realizes $T$? The most efficient general construction appears to be that of Bednay and Boz\'{o}ki (\cite{BB}, Construction 5 and Corollary 7), which shows that any tournament with $n$ vertices can be realized with a set of dice with $\lceil \frac{6}{5}n \rceil$ sides. In this paper, we give a construction which allows us to realize any tournament on $n$ vertices with dice with at most $n+1$ sides.

\section{Preliminaries} \label{sect:design}

Our construction relies on a well-known construction from graph theory, which we recall here. The edge set of the (undirected) complete graph on the vertices $\{1,2, \ldots \}$ is the set of all unordered pairs chosen from that set: $\{\{1,2\}, \{1,3\}, \ldots \{n-1,n\}\}$.

\begin{theorem} \label{thm:design} Let $E = \{\{1,2\}, \{1,3\}, \ldots \{n-1,n\}\}$. If $n$ is even, there is a partition of $E$ into $n-1$ sets of size $n/2$, where no two pairs within a single set share an element. If $n$ is odd, there is a partition of $E$ into $n$ sets of size $\frac{n-1}{2}$ where no two pairs within a single set share an element. \end{theorem}

In the language of graph theory, this theorem establishes the edge chromatic number of $K_{n}$. As such, in the following we will refer to the pairs as \textit{edges} and the elements of those pairs as \textit{vertices}. %Also note that this theorem is natural in the context of round-robin tournaments, where we have $n$ teams, each of which must play each other team once. The vertices represent teams, and the edges of the graph represent the games to be played. In the even case, we are looking for a way to schedule the games on a series of days (the sets in the partition of the edges) so that each team plays one game per day. In the odd case, such an arrangement would be impossible, so we instead look for a way to schedule the games so that each team has one day off, but plays one game each of the other days. The constructions of these partitions and the proof of their effectiveness are found in many introductory design theory texts. (See for example \cite{A}, Theorem 1.2.1. and Exercise 17 of Chapter 1.)

We do not give a full proof of the theorem, since we only need the statement of the standard construction. In the odd case, let $n=2k+1$. Then, for $i$ from 1 to $n$, we let $Y_{i}$ be the set containing the edges $\{i-1,i+1\}, \{i-2,i+2\}, \ldots \{i-k,i+k\}$. We reduce the entries of these pairs mod $n$ if needed, so that all entries are in the range from 1 to $n$. The standard proof of this theorem (see, for example \cite{A}) shows that the set $Y_{i}$ form a partition of the edges of $K_{n}$, where a vertex $j$ appears in exactly one edge in every set $Y_{i}$ except $Y_{j}$. To be able to refer to specific pairs later, we let $Y_{ij} = \{i+j,i-j\}$ for $j$ from $1$ to $k$. %Then by the same proof as above, for each $j$ from $1$ to $k$, each vertex appears exactly twice in sets of the form $Y_{ij}$. Then for each pair $(i,j)$ in $Y_{1}$, we form $Z_{i} = Y_{i} \cup Y_{j} \cup \{\{i,j\}\}$. Each of these will contain each vertex label exactly twice, since each label appears once in each of $Y_{i}$ and $Y_{j}$, except $i$ and $j$ themselves, which appear once total in $Y_{i}$ and $Y_{j}$. So adding $(i,j)$ to the set $Z_{i}$ makes them appear twice total. Then the $Z_{i}$'s will form a partition of the edges of the tournament since the $Y_{i}$ did.
As an example for the next section, note that for $n=7$, this algorithm creates the sets in Figure \ref{ex2}. %Notice now that each vertex appears twice in each column of the table.

\begin{figure}[ht] \begin{centering} \begin{tabular}{c|c} \textrm{Set:} & \textrm{Pairs:} \\ \hline $Y_{1}$ & \{2,7\},\{3,6\},\{4,5\} \\ \hline $Y_{2}$ & \{1,3\},\{4,7\},\{5,6\} \\ \hline $Y_{3}$ & \{2,4\},\{1,5\},\{6,7\} \\ \hline $Y_{4}$ & \{3,5\},\{2,6\},\{1,7\} \\ \hline $Y_{5}$ & \{4,6\},\{3,7\},\{1,2\} \\ \hline $Y_{6}$ & \{5,7\},\{1,4\},\{2,3\} \\ \hline $Y_{7}$ & \{1,6\},\{2,5\},\{3,4\} \end{tabular} \caption{The construction for $n=7$} \label{ex2} \end{centering} \end{figure}

In the even case, we let $n=2k$ and define the $i$th set of our partition, $Y_{i}$, to contain the edges $\{i,n\}, \{i-1,i+1\}, \{i-2,i+2\}, \ldots \{i-(k-1),i+(k-1)\}$, where the entries besides $n$ are reduced mod $n-1$ to fall in the range between $1$ and $n-1$. (We will not mention this reduction later, but it is always in effect when discussing these pairs.) In the case where $n=6$, this creates the partition in Figure \ref{pfex1}.

\begin{figure}[ht] \begin{centering} \begin{tabular}{c|c} \textrm{Set:} & \textrm{Pairs:} \\ \hline $Y_{1}$ & \{2,5\},\{1,6\},\{3,4\} \\ \hline $Y_{2}$ & \{1,3\},\{2,6\},\{4,5\} \\ \hline $Y_{3}$ & \{2,4\},\{3,6\},\{1,5\} \\ \hline $Y_{4}$ & \{3,5\},\{4,6\},\{1,2\} \\ \hline $Y_{5}$ & \{1,4\},\{5,6\},\{2,3\} \end{tabular}  \caption{The construction for $n=6$} \label{pfex1} \end{centering} \end{figure}

As above,  the standard proof of the theorem shows that each edge appears once in this table, and each vertex appears once in each row of this table. An additional property of these sets that we will need in the even case is that for each $j$ from $1$ to $k-1$, each vertex other than $n$ appears exactly twice in sets of the form $\{i + j, i-j\}$. To see this, note that the two vertices in $\{i+j,i-j\}$ have a difference of exactly $2j$ mod $n-1$. Thus for a vertex $m$, the sets containing $m$ are exactly $\{m,m+2j\}$ and $\{m,m-2j\}$. These sets are distinct since $m$ cannot be congruent to $m \pm 2j$ mod $n-1$ since $2j < n-1$, and if $m+2j \equiv m-2j$ mod $n-1$, then $4j \equiv 0$ mod $n-1$, and since $n-1$ is odd, this would mean $j \equiv 0$ mod $n-1$.

In the case where $n$ is even but not divisible by 4, we again define a numbering $Y_{ij}$ of the pairs in $Y_{i}$. However, we choose a somewhat non-intuitive ordering, for reasons that will be clear later. (This ordering is also reflected in the way we have arranged the pairs in Figure $\eqref{pfex1}$.) We define $Y_{ij}$ for $j$ from $1$ to $\frac{n-2}{4}$ to be the pair $\{i + j, i-j\}$. Then we let $Y_{i\frac{n+2}{4}} = \{i,n\}$. Finally, for $j$ from $\frac{n+6}{4}$ to $\frac{n}{2}$, we let $Y_{ij} = \{i+j-1,i-j+1\}$. Intuitively, we start with the pairs in the ``natural'' order: $\{i,n\}, \{i+1,i-1\}$, etc. Then we move the pair $\{i,n\}$ so it is the middle pair in our list. This creates the ordering just described.

Notice we have arranged the edges $Y_{ij}$ in order in our table, so we can restate our previous observation as saying that each vertex other than $n$ appears twice in each column other than the middle column.

There is one other property of the sets $Y_{ij}$ which will be important for our construction below. In Figure $\ref{ex2}$, note that for any two vertices $w$ and $x$, the number of rows for which $w$'s pair appears to the left of $x$'s is the same as the number of rows for which $w$'s pair appears to the right of $x$'s. For example, 3 appears to the left of 6 in rows 2 and 4, and 3 appears to the right of 6 in rows 5 and 7. (Row 1 has the pair $\{3,6\}$, and 3 and 6 are absent from rows 3 and 6,respectively, so these rows are ignored.) A similar statement is true about Figure \ref{pfex1}. We formalize these observations in the following lemma.

\begin{lemma} \label{lemma:lemma} If $n$ is even and not a multiple of 4, let $w$ and $x$ be two vertices. For $i$ from 1 to $n$, define $w_{i}$ and $x_{i}$ so that $w \in Y_{iw_{i}}$ and $x \in Y_{ix_{i}}$. Then $w_{i} < x_{i}$ for exactly $\frac{n}{2}-2$ values of $i$.

If $n$ is odd, let $w$ and $x$ be any two vertices. For $i$ from 1 to $n$, define $w_{i}$ and $x_{i}$ so that $w \in Y_{iw_{i}}$ and $x \in Y_{ix_{i}}$, assuming such values exist. Then $w_{i} < x_{i}$ for exactly $\frac{n-3}{2}$ values of $i$.
\end{lemma}

\begin{figure}[ht]
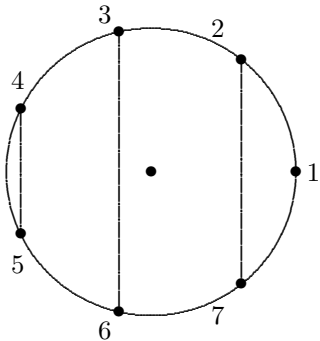

\[ \beginpicture
\setcoordinatesystem units <0.15in,0.15in>
\circulararc 360 degrees from 5 0 center at 0 0
\put{$\bullet$} at 0 0
\put{$\bullet$} at 5 0
\put{1} at 5.6 0
\put{$\bullet$} at 3.11 3.90
\put{2} at 2.3 5
\put{$\bullet$} at 3.11 -3.90
\put{7} at 2.3 -5
\plot 3.11 3.90 3.11 -3.90 /
\put{$\bullet$} at -1.11 4.874
\put{3} at -1.6 5.5
\put{$\bullet$} at -1.11 -4.874
\put{6} at -1.6 -5.5
\plot -1.11 4.874 -1.11 -4.874 /
\put{$\bullet$} at -4.5 2.169
\put{4} at -4.6 3.2
\put{$\bullet$} at -4.5 -2.169
\put{5} at -4.6 -3.2
\plot -4.5 2.169 -4.5 -2.169 /
\setlinear
\endpicture \]
\caption{The set $Y_{1}$ in the case $n=7$}
\label{pfex3}
\end{figure}

Proof. We begin with the odd case. Imagine a circle with the $n$ points $1,2, \ldots n$ drawn, in order, equally spaced around the circle. Then the pairs $Y_{ij}$ are, in order from $j=1$ to $j=k$, the two vertices closest to $i$, then the next two vertices, etc. (See Figure \ref{pfex3}.)

\begin{figure}[ht]
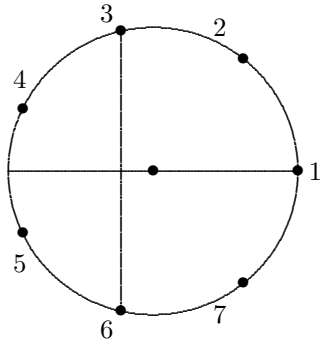

\[ \beginpicture
\setcoordinatesystem units <0.15in,0.15in>
\circulararc 360 degrees from 5 0 center at 0 0
\put{$\bullet$} at 0 0
\put{$\bullet$} at 5 0
\put{1} at 5.6 0
\plot 5 0 -5 0 /
\put{$\bullet$} at 3.11 3.90
\put{2} at 2.3 5
\put{$\bullet$} at 3.11 -3.90
\put{7} at 2.3 -5
\put{$\bullet$} at -1.11 4.874
\put{3} at -1.6 5.5
\put{$\bullet$} at -1.11 -4.874
\put{6} at -1.6 -5.5
\plot -1.11 4.874 -1.11 -4.874 /
\put{$\bullet$} at -4.5 2.169
\put{4} at -4.6 3.2
\put{$\bullet$} at -4.5 -2.169
\put{5} at -4.6 -3.2
\setlinear
\endpicture \]
\caption{Comparing 3 and 6 in the case $n=7$}
\label{pfex4}
\end{figure}

So for a vertex $w$, the number $w_{i}$ is exactly the ``distance'' from $w$ to $i$ on this circle. (We count moving from one vertex to the next as a distance of 1.) Thus for a pair of vertices $w$ and $x$, $w_{i} < x_{i}$ exactly if $w$ is closer to $i$ than $x$ is. The pair $\{w,x\}$ will appear in the set $Y_{m}$ for the unique $m$ which is equidistant from both of them. In geometric terms, the perpendicular bisector of the segment between $w$ and $x$ will be the diameter of the circle going through $m$, and this line has half of the remaining points on either side of it. Thus, $w$ is closer to $i$ than $y$ is for the $\frac{n-3}{2}$  other points on $w$'s side of the bisector. Thus $w_{i} < y_{i}$ for exactly $\frac{n-3}{2}$ values of $i$. In Figure \ref{pfex4}, we see that $\{3,6\}$ is in $Y_{i}$ in Figure \ref{ex2} for $n=7$. Also, as observed, 3 appears to the left of 6 in the second and fourth rows, and to the right of 6 in the fifth and seventh rows of Figure  \ref{ex2}.

\begin{figure}[ht]
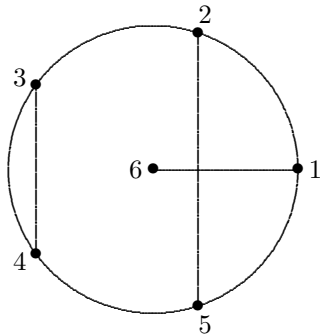

\[ \beginpicture
\setcoordinatesystem units <0.15in,0.15in>
\circulararc 360 degrees from 5 0 center at 0 0
\put{$\bullet$} at 0 0
\put{6} at -.6 0
\put{$\bullet$} at 5 0
\put{1} at 5.6 0
\put{$\bullet$} at 1.545 4.75
\put{2} at 1.8 5.4
\put{$\bullet$} at 1.545 -4.75
\put{5} at 1.8 -5.4
\put{$\bullet$} at -4.045 2.93
\put{3} at -4.6 3.2
\put{$\bullet$} at -4.045 -2.93
\put{4} at -4.6 -3.2
\setlinear
\plot 5 0 0 0 /
\plot 1.545 4.75 1.545 -4.75 /
\plot -4.045 2.93 -4.045 -2.93 /
\endpicture \]
\caption{The set $Y_{1}$ in the case $n=6$}
\label{pfex2}
\end{figure}

For the even case of our proof, we can draw a similar picture, now putting vertex $n$ as the center of the circle and spacing the points $1,2, \ldots n-1$ evenly on the circle. Then $Y_{i}$ can be visualized by connecting $i$ to $n$ and connecting the other points pairwise with lines perpendicular to the segment between $i$ and $n$, as in Figure \ref{pfex2}.

If the edges $Y_{ij}$ were in the ``natural'' order in this case, an argument analoguous to the odd case would suffice to prove the lemma for $w$, $x \neq n$. However, the pair $\{i,n\}$ has been moved to be the ``middle'' pair $Y_{i\frac{n+2}{4}}$ for each $i$. So we just need to check that  the result holds when $w$ or $x$ is $n$, and that moving the pairs $\{i,n\}$ does not change the result when $w,x \neq n$. The first observation follows from the fact that, for fixed $j \neq \frac{n+2}{4}$, any vertex $w$ occurs in exactly two sets $Y_{ij}$, as observed above. For the second, notice that for any vertices $w$ and $x$ other than $n$, moving the pair $\{i,n\}$ does not interfere with the relative order of the pairs containing $w$ and $x$ in $Y_{i}$ unless $i=w$ or $x$. However, the pair containing $w$ in $Y_{x}$ is the pair $\{x+j,x-j\}$, where $w \cong x \pm j$ mod $n-1$. But then of course $x \cong w \pm j$ mod $n-1$, so that $w_{x} = x_{w}$. And since $w_{w} = x_{x} = \frac{n+2}{4}$, we know that $w_{x} < x_{x}$ exactly if $x_{w} < w_{w}$. Thus moving the pair $\{i,n\}$ as described will change the relative orders of $w$ and $x$ either in both of rows $w$ and $x$ or neither, and the result holds for our final ordering.

$\square$

Note: The ordering of the $Y_{ij}$ in the even case was chosen exactly so that the lemma holds in the case where $w$ or $x$ is $n$. This will be important for our main construction.

\section{The Construction}

We need one piece of terminology for our main proof. For two dice $X_{1}$ and $X_{2}$, the number of \textit{face wins} (or simply \textit{wins}) for $X_{1}$ over $X_{2}$ is the number of ordered pairs $(a,b)$ where $a$ is a face of $X_{1}$, $b$ is a face of $X_{2}$, and $a > b$. So in our example $$\begin{tabular}{c|c} \textrm{Die:} & \textrm{Faces:} \\ \hline $X_{1}$ & 1,5,9 \\ $X_{2}$ & 3,4,8 \\ $X_{3}$ & 2,6,7 \end{tabular},$$ $X_{1}$ has 5 face wins over $X_{2}$, corresponding to the pairs $(9,8),(9,4),(9,3),(5,4)$, and $(5,3)$. Note that if the dice have $n$ sides, then the probability that $X_{1}$ rolls higher than $X_{2}$ is exactly the number of wins for $X_{1}$ over $X_{2}$ divided by $n^{2}$.

We are now in a position where we can prove our main theorem. %We construct our dice using the columned construction, and each column of the table we construct will correspond to one set $Y_{i}$, so that that particular face of the corresponding dice is what makes the matchups between pairs in $Y_{i}$ have the desired winner.
We construct our set of $n$ dice so that the first face of each die has a number from 1 to $n$, the second face has a number from $n+1$ to $2n$, etc. Moreover, the $i$th face of each die will be constructed based on the set $Y_{i}$, and for each pair of dice in $Y_{i}$, the matchup between those two dice is determined by that particular face. So, for example, when $n=6$, $Y_{1} = \{\{1,6\},\{2,5\},\{3,4\}\}$. This means that we construct our dice so that $X_{1}$'s first face is greater than $X_{6}$'s first face exactly if 1 beats 6 in our tournament $T$. The same holds for $X_{2}$ and $X_{5}$, and $X_{3}$ and $X_{4}$. The difficulty lies in ensuring that the remaining sides give a probability of exactly $1/2$ that $X_{1}$ beats $X_{6}$, so that the first face determines the matchup between the dice as intended. It turns out that the partitions above give us a way to do just that.

\begin{theorem} \label{thm:main} Let $T$ be a tournament on $n$ vertices. If $n$ is odd, there is a set of $n$-sided dice that realize $T$. If $n$ is divisible by 4, there is a set of $n+1$-sided dice that realize $T$. If $n$ is even and not divisible by 4, then there is a set of $n-1$-sided dice that realize $T$.  \end{theorem}

Proof. First assume $n$ is odd.  Construct the sets $Y_{i}$ as in Theorem $\ref{thm:design}$. Now we construct our dice as follows. For $i$ from 1 to $n$, the $i$th face of each die will contain a number from $n(i-1)+1$ to $ni$. Specifically, if $i=k$, the $i$th face of Die $k$ will contain $n(i-1) +1$. If $i \neq k$, then there is some set $Y_{ij}$ containing $k$. In this case, the $ith$ face of Die $k$ will contain either $n(i-1) + 2j$ or $n(i-1) + 2j+1$. Then for each pair $Y_{ij}$, the $i$th face of the corresponding dice will have one of the two numbers $n(i-1) + 2j$ or $n(i-1) + 2j+1$, and we give the higher number to the die designated to win the matchup by the corresponding edge in $T$.

As an example of this construction, Figure \ref{pfex5} is the constructed set of dice for $n=7$. In this table, the $i$th column represents the $i$th face of each die. As a shorthand, the numbers in the $i$th column are reduced mod 7, but they represent the face labels $7(i-1) + 1$ through $7i$. The notation $x/y$ means that the particular label is chosen from between those two values according to $T$ as described above.

\begin{figure}[ht] \begin{centering}\begin{tabular}{|c|c|c|c|c|c|c|c|}\hline Die 1: & 1 & 2/3 & 4/5 & 6/7 & 6/7 & 4/5 & 2/3 \\ \hline Die 2: &   2/3&1  & 2/3 & 4/5 & 6/7 & 6/7 & 4/5 \\ \hline Die 3: & 4/5 & 2/3 &1 & 2/3 & 4/5 & 6/7 & 6/7 \\ \hline Die 4: &  6/7& 4/5& 2/3 &1 & 2/3 & 4/5 & 6/7 \\ \hline Die 5: & 6/7& 6/7 & 4/5 & 2/3 & 1& 2/3 & 4/5 \\ \hline Die 6: & 4/5 & 6/7 & 6/7 & 4/5 & 2/3 &1 & 2/3 \\ \hline Die 7: & 2/3 & 4/5 & 6/7 & 6/7 & 4/5 & 2/3 & 1\\ \hline  \end{tabular}  \caption{The constructed table for $n=7$} \label{pfex5} \end{centering}\end{figure}

By construction, each die gets at least $\binom{n}{2}$ wins over each other die, since the $i$th face of any die will always beat the $k$th face of any die when $k < i$. Also, by the properties of the constructed partitions, a die $w$ is guaranteed an extra $\frac{n-3}{2}$ wins over each other die $x$ from the $\frac{n-3}{2}$ values of $i$ where $w_{i} > x_{i}$ (as defined in Lemma \ref{lemma:lemma}), since that guarantees that Die $w$ will have a higher number in the corresponding column of this table. Also, Die $w$ is guaranteed to beat Die $x$ in column $x$, since Die $x$ has a 1 there. This is a total of $\frac{n^{2}-1}{2}$ wins for each die against each other die. Then if $\{w,x\} = Y_{ij}$, then the overall matchup between Die $w$ and Die $x$ is determined by which die has the higher number on the $i$th face. Since we chose that face to match $T$, this set of dice will realize $T$.

For example, consider the tournament on 7 vertices where die $X_{i}$ beats die $X_{j}$ whenever $i<j$, except that $X_{7}$ beats $X_{1}$. (This tournament is ``almost transitive'' - if $X_{1}$ beat $X_{7}$, it would be transitive.) Thus, in almost every column, for every pair of undetermined entries in the table in Figure \ref{pfex5}, the lower-numbered die gets the higher value, so that the lower-numbered die beats the higher-numbered one. The exception to this is the fourth column, where we choose 6 for row 1 and 7 for row 7, so that $X_{7}$ beats $X_{1}$. The resulting set of dice is shown in Figure $\ref{pfex6}$:

\begin{figure}[ht] \begin{centering}\begin{tabular}{|c||c|c|c|c|c|c|c|}\hline $X_{1}$: & 1 & 10 & 19 & 27 & 35 & 40 & 45 \\ \hline $X_{2}$: & 3 & 8 & 17  & 26 & 34 & 42 & 47 \\ \hline $X_{3}$: & 5 & 9 & 15 & 24 & 33 & 41 & 49 \\ \hline $X_{4}$: & 7 & 12 & 16 & 22 & 31 & 39 & 48 \\ \hline $X_{5}$: & 6 & 14 & 18 & 23 & 29 & 38 & 46 \\ \hline $X_{6}$: & 4 & 13 & 21 & 25 & 30 & 36 & 44 \\ \hline $X_{7}$: & 2 & 11 & 20 & 28 & 32 & 37 & 43\\ \hline  \end{tabular}  \caption{A set of dice realizing the almost transitive tournament for $n=7$} \label{pfex6} \end{centering}\end{figure}

The remaining cases are slight variations on this procedure. If $n$ is divisible by 4, we make a new tournament $T'$ by adding a die which beats all other dice. Then we apply the above algorithm to construct a set of $n+1$ columned $n+1$-sided dice that realize $T'$. By deleting the added die from this set, we can construct a set of $n$ columned $n+1$-sided dice that realize $T$.

We could do the same if $n$ is even but not divisible by 4, but it is possible in this case to construct a set of dice with $n-1$ sides that realize $T$. %We will use the columned construction as before. We make the table of pairs of vertices as in \eqref{table}, where now the entry in row $j$ and column $i$ (counting from row 0) if the set $Y_{ij}$ from the even case of Theorem $\ref{thm:design}$. However, we then move row 0 to the middle of the table, so it is now row $\frac{n-2}{4}$, sliding the intervening rows up:
To do so, we use the even case of the algorithm of Theorem \ref{thm:design} to construct the sets $Y_{ij}$. As before, for $i$ from 1 to $n-1$, the $i$th face of each die will contain a number from $n(i-1)+1$ to $ni$. If $k \in Y_{ij}$, then the $i$th face of Die $k$ will contain either $n(i-1) + 2j-1$ or $n(i-1) + 2j$.  Then as before, for each pair $Y_{ij}$, the $i$th face of the corresponding dice will have one of two numbers, and we give the higher number to the die designated to win the matchup by the corresponding edge in $T$.The result for $n=6$ is shown, in Figure \ref{pfex7}, using the same notation as Figure $\ref{pfex5}$.

\begin{figure}[ht] \begin{centering}\begin{tabular}{|c|c|c|c|c|c|} \hline Die 1: & 3/4 & 1/2 & 5/6 & 5/6 & 1/2 \\ \hline Die 2: & 1/2 & 3/4 & 1/2 & 5/6 & 5/6 \\ \hline Die 3:  & 5/6 & 1/2 & 3/4 & 1/2 & 5/6 \\ \hline Die 4: & 5/6 & 5/6 & 1/2 & 3/4 & 1/2 \\ \hline Die 5:  & 1/2 & 5/6 & 5/6 & 1/2 & 3/4 \\ \hline Die 6:  & 3/4 & 3/4 & 3/4 & 3/4 & 3/4  \\ \hline  \end{tabular}  \caption{The constructed table for $n=6$} \label{pfex7} \end{centering}\end{figure}

Then by the exact same logic as before (using Lemma \ref{lemma:lemma} again) we can count that every die is guaranteed $\frac{(n-1)^{2}-1}{2}$ wins against every other die, and the matchup between any two dice in a pair $Y_{ij}$ is determined by which die gets the higher number in column $j$, which we chose to match $T$. $\square$

Note also that the sets of dice constructed by this process are uniform in the sense that if $X_{i} \succeq X_{j}$, then the probability that die $X_{i}$ rolls higher than $X_{j}$ is exactly $\frac{1}{2} + \frac{1}{2k^{2}}$, where $k$ is the number of sides on the dice. This is the probability closest to 1/2 that is achievable on $k$-sided dice, and every matchup is decided with this probability. (This is a slightly stronger condition than the notion of ``balanced'' introduced in \cite{SS}).

\section{Conclusion}

There have been a number of general algorithms discovered that create a set of dice that realize a given tournament. The previous algorithm with the fewest guaranteed number of sides (see \cite{BB}) was inherently inductive, adding more and more dice until the desired tournament is achieved. Our algorithm has the advantages of being direct and slightly more efficient. There is also an algorithm described by Schaefer (\cite{S}) for starting with any set of dice, and adding faces to the dice to change the edges until the desired tournament is achieved. Naively, starting from a transitive tournament realized by 1-sided dice, it is always possible to realize any tournament on $n$ vertices with at most $2n-1$ sides on the dice using Schaefer's algorithm. But it seems likely that for a given tournament, there are more intelligent ways to use the algorithm to achieve that tournament in fewer sides. (The general problem of finding how many edges must be changed to turn a given tournament into a transitive tournament is known as the \textit{feedback arc set problem}, and is, in general,  quite difficult.) Also, if the realization of Paley tournaments on $p$ in \cite{BB} using dice with $\frac{p-1}{2}$ sides is any indication, our upper bound of $n$ sides is probably rarely, if ever, the best possible.  Thus the question of the smallest size of dice needed to realize a given tournament is ripe for further study.

\end{document}